
\documentstyle{amsppt}

\def\O{\Cal O}
\def\g{\goth g}
\def\lra{\longrightarrow}
\def\End{\Cal{E}\text{\rm nd}}
\def\tr{\text{\rm tr}}
\def\im{\text{\rm im}}
\def\ker{\text{\rm ker}}
\def\X#1{X\kern-.5em<\kern-.3em#1\kern-.3em>}

\magnification=1200

\topmatter

\title
On The Local Geometry of Moduli Spaces of Vector Bundles
\endtitle

\author
Ziv Ran
\endauthor

\email
ziv\@math.ucr.edu
\endemail

\date
April 24, 1995
\enddate

\endtopmatter

\document

In [R1] the author gave a canonical construction of the universal deformation
of a compact complex manifold.  Here we shall give an analogous
construction for the universal deformation of a simple vector bundle
(over a fixed base); this construction differs somewhat from that of [R1],
featuring a systematic use of the Jacobi complexes of a module.
The resulting construction of universal or Poincar\'e bundle is
applicable in other problems as well, such as deformations of manifold,
where it yields a simplification of the method employed in [R1] to construct
the structure sheaf of the total space of the universal deformation.

As an application of our method, we shall prove the closedness of
{\it trace forms}---these are $H^2(\O_X)$--valued 2--forms on the
moduli space of bundles, $E$, on $X$ induced by the trace or Killing form
on $\End(E)$---
$$
\tau:  H^1(\End(E))\otimes H^1(\End(E)) \rightarrow
H^2(\End(E)) {\buildrel\tr\over\rightarrow} H^2(\O_X).
$$

For $X$ a surface with trivial canonical bundle, $\tau$ is in fact,
by Serre duality, a non--degenerate scalar valued 2--form, hence,
being closed, yields a symplectic structure on the moduli space, a result
first proven by Mukai [M];  see also Kobayashi's book [K] for an analytic
approach.

\medskip

\subheading{Acknowledgment}
I am grateful to Professors Maruyama and Mukai for their invitation to attend
the symposium and to Professor Fujiki for suggesting the problem of
closedness of trace forms.

\medskip

\subheading{1.  Calculus}

\medskip

\subsubhead
{1.1.  Enveloping Algebras and Differential Operators}
\endsubsubhead

\medskip

Let $\g$ be a Lie algebra, say over a field
$\Bbb C$ of characteristic 0, and $E$ a faithful $\g$--module.
As is well known, the universal enveloping algebra $\Cal U=\Cal U(\g)$
may be identified as vector space, though not as algebra, with
$\oplus^\infty_{i=0}S^i(\g)$, and this identification takes
$\Cal U^m \subset \Cal U$, the elements of order $\le m$, to
$\oplus^m_{i=0}S^i(\g)$.  The algebra structure on $\Cal U(\g)$
may be described as follows.
Let $\Cal S^\cdot = \Cal S^\cdot(E) = \oplus^\infty_{i=0}S^i(E)$ be the
symmetric algebra on $E$ and note that elements of $\g$ acting on $E$
extend to graded (i.e. homogeneous of degree zero) derivations on
$\Cal S^\cdot$; in fact it is easy to see that the Lie algebra of
all graded derivations coincides with $\goth{gl}(E) = \text{\rm End}(E)$.
Therefore if we denote by $\Cal D$ (resp. $\Cal D^m$) the algebra
(resp. $\g$--module) of all (resp. all degree $\le m$) graded
$\Bbb C$--linear differential operators on $\Cal S^\cdot$, we
obtain an injective filtration preserving algebra and $\g$--module
homomorphism
$$
\rho:\Cal U(\g) \rightarrow \Cal D.
$$
Moreover it is easy to see that for $\g = \goth{gl}(E)$, $\rho$ is
in fact an isomorphism.

We want now to identify $\im(\rho)$ in case
$\g=\goth{sl}(E)$, the traceless endomorphisms (assuming $E$ is
finite--dimensional).
To this end, note that a graded endomorphism $\varphi^\cdot$ of $\Cal S^\cdot$
admits a graded trace $\tr^\cdot(\varphi^\cdot)$, which itself is a
graded endomorphism of $\Cal S^\cdot$ given by
$$
\multline
\tr^i(\varphi^\cdot) =
\text{\rm image of~}
(\varphi^{i+1} \in S^{i+1}(E)\otimes S^{i+1}(E^*)\subset
S^i{E}\otimes S^i(E^*)\otimes E\otimes E^*)
\\
\text{\rm ~in~}
S^i(E)\otimes S^i(E^*)
\text{\rm ~under~} \text{\rm id}\otimes\text{\rm id}\otimes\tr
\endmultline
$$
In other words, $\tr^\cdot$ is the natural extension of
$\tr: E\otimes E^* \rightarrow \Bbb C$
as a derivation of degree $-1$ on the algebra of all graded
endomorphisms of $\Cal S^\cdot$.
Clearly $\tr^\cdot(\varphi^\cdot) = 0$ iff for all
$e_1,\ldots,e_i \in E$, $e_1^*,\ldots,e_i^* \in E^*$, acting on
$\Cal S^\cdot$ in the standard manner (interior multiplication), the
trace (in the usual sense) of
$e_1^*\ldots e_i^*\circ \varphi^{i+1}\circ e_1 \ldots e_i$ is zero.
 From this it is easy to see that $\rho(\Cal U(\goth{sl}(E)))$ coincides
with the sub--algebra $\Cal {SD} \subset \Cal D$ of traceless operators.

The above considerations extend immediately to the case where $E$ is a
vector bundle over a $\Bbb C$--ringed space $(X,\O_X)$ and $\g$
is a Lie sub--algebra of $\goth {gl}(E) = \End_{\O_X}(E)$.

\medskip

\subsubhead
{1.2.  Functors on $S$--Modules}
\endsubsubhead

\medskip

In [R1] we showed how an Artin local $\Bbb C$--algebra may be reconstructed
from a certain `order--symbolic' or OS structure on the space of
$\Bbb C$--valued differential operators on $S$.
Our purpose here is to note an analogue of this for $S$--modules.

Now fix a local $\Bbb C$--algebra $S$ with maximal ideal $\goth m$ and
residue field $S/\goth m = \Bbb C$ and put
$$
B_0^i = D^i(S,\Bbb C) = \text{\rm Hom}(S_i,\Bbb C) = S_i^*,
\kern2em S_i = S/\goth m^{i+1},
$$
and
$$
B^i = D_+^i(S,\Bbb C) = (\goth m/\goth m^{i+1})^*.
$$
For an $S$--module $E$, put
$$
B^i(E) = B_0^i\otimes_S E,
$$
where $B^i_0$ is viewed as $S$--bi--module and $E$ as (symmetric)
$S$--bi--module.  At least when $E\otimes S_i$ is $S_i$--free,
$B^i(E)$ may be identified with the right $S$--module of differential
operators $D^i(E^\vee,\Bbb C)$, $E^\vee=\text{\rm Hom}_S(E,S)$.

We have a symbol map
$$
\sigma^i: B_0^i \rightarrow B^i\otimes_{\Bbb C}B^{i-1}_0
$$
which factors through $F_i(B^i\otimes_{\Bbb C}B^{i-1})$, where
$F_i$ is the filtration induced by the order filtration on $B^i$,
and this gives rise to a symbol map
$$
\sigma^i_E: B^i(E) \rightarrow B^i\otimes_{\Bbb C}B^{i-1}(E),
$$
which again factors through $F_i(B^i\otimes_{\Bbb C}B^{i-1}(E))$.
These $\sigma^i_E$, $i\le m$, together with the obvious maps
$B^0(E) \rightarrow B^1(E)\rightarrow \ldots\rightarrow B^m(E)$
are referred to as a `modular order--symbolic' (MOS) structure on
$B^m(E)$.  Note that $B^m(E)$ itself is a right $S_m$--module,
called the $m$--th transpose of $E$.

``Dually,'' suppose we are given an MOS structure $G^\cdot$, $G^i$ a right
$S_i$--module.  We then define an $S_m$--module $C^m(G^\cdot)$, called
the module of quasi--scalar homomorphisms $B^m_0 \rightarrow C^m(G^\cdot)$,
inductively as follows.
$$
C^0(G) = G^0
$$
$$
C^i(G) = \text{\rm all right~} S\text{\rm --linear maps~}
\varphi^i: B^i_0 \rightarrow G^i
$$
such that for some $\varphi^{i-1}: B^{i-1}_0\rightarrow G^{i-1}$ the following
diagrams commute.
$$
\CD
B_0^{i-1} @>{\varphi^{i-1}}>> G^{i-1} &&&&
\kern10em&& B_0^i @>{\varphi^i}>> G^i \\
@VVV @VVV &&
\kern10em&& @V{\sigma^i}VV @VV{\sigma^i_G}V  \\
B_0^i @>{\varphi^i}>> G^i &&&&
\kern10em&& B^i\otimes B^{i-1}_0 @>{\text{\rm id}\otimes\varphi^{i-1}}>>
B^i\otimes G^{i-1} \\
\endCD
$$

Note that we have natural maps
$$
E\rightarrow C^m(B^m(E))
$$
$$
B^m(C^m(G^\cdot)) \rightarrow G^\cdot.
$$
At least when $E$ is $S_m$--free (resp. $G^\cdot$ is `co--free',
i.e. a sum of copies of $B^m_0$ with the standard MOS structure),
these are isomorphisms.

\medskip

\subsubhead
{1.3.  Exterior Derivative}
\endsubsubhead

\medskip

Our purpose here is to give a convenient interpretation of the Cartan
formula for exterior derivative.  Let $M$ be a manifold with tangent
sheaf $T$, and consider the sheaf $B^1(\bigwedge^i T)$,
which may be identified as the sheaf of $\O_M$--valued
first--order differential operators on the dual
$\bigwedge^i T^\vee = \Omega^i$, and fits in an exact sequence
$$
0\rightarrow \bigwedge^i T\rightarrow B^1(\bigwedge^i T)
\rightarrow  T\otimes\bigwedge^i T\rightarrow 0.
$$
Over the sub--sheaf $\bigwedge^{i+1} T\subset T\otimes\bigwedge^i T$,
this sequence admits a canonical splitting
$\varphi_i:\bigwedge^{i+1} T\rightarrow B^1(\bigwedge^i T)$
given by
$$
\multline
\varphi_i(v_0\wedge\ldots\wedge v_i)(\omega) =
\sum (-1)^jv_j\omega(v_0,\ldots,\hat v_j,\ldots,v_i)+
\\
\sum_{j<k}(-1)^{j+k}\omega([v_j,v_k],v_0,\ldots,\hat v_j,\ldots,\hat v_k,
\ldots,v_i)
\endmultline
\tag{*}
$$
for $\omega \in \Omega^i$.  The Cartan formula says in particular
that $\varphi_i$ is $\O_M$--linear.  Now given an $i$--form, $\omega$,
viewed as a map $\bigwedge^i T \rightarrow \O_M$, it extends
to a map $B^1(\bigwedge^i T)\rightarrow B^1(\O_M) = \O_M\oplus T$.
Projecting on the $\O_M$ factor, we get a map
$$
\tilde\omega: B^1(\bigwedge^i T)\rightarrow \O_M.
$$
Cartan's formula then says that the exterior derivative $d\omega$
coincides with $\tilde\omega\circ\varphi_i$.

Note that $\im(\varphi_1)$ coincides with the kernel of the
natural map
$$
{}^td_0: B^1( T)\rightarrow  T^2=D^2_+(\O,\O)
$$
$$
{}^td_0(a) = a\circ d_0
$$
where $a:\Omega^1\rightarrow\O$ is a first--order differential operator and
$d_0:\O\rightarrow \Omega^1$ is exterior derivative.  Indeed we have an
exact sequence
$$
\matrix
&&&&&&\bigwedge^2 T \\
\\
&&&&& {\buildrel\varphi_1\over\swarrow} &\downarrow \\
\\
0&\rightarrow& T&\rightarrow&B^1( T)&\rightarrow& T\otimes T&
\rightarrow&0 \\
\\
&&\Vert&&{}^td_0\downarrow&&\downarrow \\
\\
0&\rightarrow& T&\rightarrow& T^2&\rightarrow&S^2 T&
\rightarrow&0 \\
\endmatrix
$$
Similarly, using the inclusion
$B^1(\bigwedge^i T)\subset
B^1( T\otimes\bigwedge^{i-1} T) =
B^1( T)\otimes_\O\bigwedge^{i-1} T$, we may identify
$\im(\varphi_i)$ as the kernel of the natural map
$B^1(\bigwedge^i T){\buildrel {}^td_0\otimes\text{\rm id} \over \lra}
 T^2\otimes\bigwedge^{i-1} T = B^2(\bigwedge^{i-1}T)/\Omega^{i-1}
=: B_+^2(\bigwedge^{i-1}T)$,
which in fact coincides with the composite
$B^1(\bigwedge^i T){\buildrel {}^td_{i-1} \over \lra}
B^2(\bigwedge^{i-1}T) \rightarrow B_+^2(\bigwedge^{i-1}T)$.

\medskip

\subheading{2.  Jacobi Complexes and Universal Deformations}

\medskip

Fix a base space $X$ which for convenience we assume to be a compact
complex space (although the construction works more generally), and a
simple vector bundle $E$ on $X$.  Let $\g=\goth{sl}(E)$, the
Lie algebra of traceless endomorphisms of $E$, so that $H^0(\g)=0$.
As in [R1], we have Jacobi complexes $J^\cdot_m(\g)$ which may be
described as follows.  Let $\X m$ be the $m$--fold very
symmetric product of $X$, i.e. the space of subsets of $X$ of
cardinality $\in [1,m]$, with the topology induced by the natural map
$X^m \rightarrow \X m$.  For $i\le m$ let $\lambda^i(\g)$ be
the image of the exterior alternating product of $\g$, supported on
$\X i \subset \X m$.  Then the bracket on $\g$ gives
rise to a map $\lambda^i(\g)\rightarrow \lambda^{i-1}(\g)$,
and these fit together to form the complex $J^\cdot_m(\g)$---
$$
\lambda^m(\g)\rightarrow\lambda^{m-1}(\g)\rightarrow
\ldots\rightarrow\lambda^1(\g) = \g,
$$
in which we put $\lambda^i(\g)$ in degree $-i$.  Similarly, the
action of $\g$ on $E$ gives rise to a complex
$J^\cdot_m(\g,E)$ on $\X m\times X$---
$$
\lambda^m(\g)\boxtimes E\rightarrow\ldots\rightarrow
\g\boxtimes E \rightarrow E,
$$
in degrees $\in [-m,0]$, where the last term $E$ is supported on the
diagonal in $\X1\times X = X\times X$.  The natural maps
$$
J_i^\cdot(\g)\rightarrow J^\cdot_m(\g),\kern5em i\le m,
$$
$$
J^\cdot_m(\g)\rightarrow F_m(J^\cdot_{m-1}(\g)\boxtimes
J^\cdot_{m-1}(\g))
$$
give rise to an OS structure on $V_m = \Bbb H^0(J^\cdot_m(\g))$ which,
as in [R1], yields a $\Bbb C$--algebra structure on
$R_m=\Bbb C\oplus V_m^*$, which turns out to be the ring of the universal
$m$--th order deformation of $E$ (see below).  To get the universal
or Poincar\'e bundle on $X\times\text{\rm Spec}R_m$, consider
$$
G^m=\Bbb R^0p_{2*}(J^\cdot_m(\g,E)).
$$
This forms a sheaf of $\O_X\otimes R_m$--modules, and the natural map
$$
J^\cdot_m(\g,E)\rightarrow F_m(J^\cdot_m(\g)\boxtimes
J^\cdot_{m-1}(\g,E))
$$
endows $G^m$ with a MOS structure compatible with the OS structure
on $V^m$, whence as in $\S$ 1.2 a sheaf
$$
\Cal P_m(E) := C^m(G^m).
$$
As $G^m$ is locally isomorphic to $(V_m\oplus\Bbb C)\otimes \O_X$, clearly
$\Cal P_m(E)$ is a locally free $R_m\otimes \O_X$ module, and it turns out
to be the Poincar\'e bundle for $E$.  (The last construction may be
applied in other situations, such as deformations of manifolds, where
it yields a simplification of the construction in [R1] of the structure
sheaf of the universal deformation of a manifold $X$.)

Notwithstanding the simple construction of $\Cal P_m(E)$ and $R_m$, the
proof of their universality involves the somewhat complicated construction
of Kodaira--Spencer bi--complexes associated to a given deformation of
$E$.  Given such a deformation, i.e. a locally free sheaf $\Cal E$ on
$X\times \text{\rm Spec}(S)$, $S$ Artin local, with
$\Cal E\otimes k(0) {\buildrel \sim\over\rightarrow} E$, we
construct bi--complexes $W_m^{\cdot,\cdot}$, and $K_m^{\cdot,\cdot}$,
which fit in exact (up to quasi--isomorphism) sequences
(in which $B^m = D^m_+(S,\Bbb C)$)---
$$
0\rightarrow J^\cdot_m(\g)\rightarrow
W_m^{\cdot,\cdot}\rightarrow B^m[1-m,m]\rightarrow 0
\kern3em\text{\rm on~}\X m,
\tag{2.1}
$$
$$
0\rightarrow J^\cdot_m(\g,E)\rightarrow
K_m^{\cdot,\cdot}\rightarrow B_S^m(\Cal E)[1-m,m]\rightarrow 0
\kern3em\text{\rm on~}\X m\times X,
\tag{2.2}
$$
whose associated coboundary maps
$$
\alpha_m:B^m\rightarrow \Bbb H^0(J^\cdot_m(\g))=V_m,
$$
$$
\beta_m:B^m_S(\Cal E)\rightarrow\Bbb R^0p_{2*}(J^\cdot_m(\g,E))=G_m
$$
respect the respective OS structures, hence give rise to a ring homomorphism
$t_{\alpha_m}:R_m\rightarrow S_m=S/\goth m^m$ and
to an $S_m$--linear map
$t_{\beta_m}:\Cal E_m:=\Cal E\otimes S_m\rightarrow
\Cal P_m(E)\otimes_{R_m}S_m$,
which turns out to be an isomorphism.

We sketch the construction of these bi--complexes.
Let $\Cal S^\cdot(E)$ be the symmetric algebra on $E$ and similarly
for $\Cal E$ and $\Cal E_m$ (as $\O_X\otimes S$--modules).
A graded, $\O_X$--linear map
$\varphi:\Cal S^\cdot(\Cal E_i)\rightarrow \Cal S^\cdot(\Cal E_j)$
is said to be a differential operator of degree $\le k$ if
for all $e^*\in E^*\subset \Cal E^*_i$ and $\goth m_{e^*}$, the ideal
in $\Cal S^\cdot(\Cal E_i)$ generated by
$\ker(e^*) \subset \Cal E_i$ and $\goth m$,
we have $\varphi(\goth m^l_{e^*})\subset \goth m^{l-k}_{e^*}$; similarly for
$\psi:S_i\rightarrow S_j$.  $\varphi$ is said to be traceless if it
can be locally written as $\sum \lambda_r\psi_r$, with
$\lambda_r:\Cal S^\cdot(E)\rightarrow S^\cdot(E)$ traceless,
$\psi_r:S_i\rightarrow S_j$.
Define
$$
\multline
\Cal A^{i,j}=\Cal A^{i,j}(\Cal E)=\lbrace(\varphi,\psi)|
\varphi:\Cal S^\cdot(\Cal E_i)\rightarrow S^\cdot(\Cal E_j),
\psi:S_i\rightarrow S_j,
\\
\O_X\text{\rm --linear graded traceless (resp.~}\Bbb C
\text{\rm --linear) diff. op. of order~}\le i-j
\\
\text{\rm such that~}\varphi|_{S_i}=\psi,\text{\rm and~} \psi(1)=0\rbrace.
\endmultline
$$
Also define
$$
\Cal A_0^{i,j}=A^{i,j}\oplus \Bbb C.
$$
These are $(S_j,S_i)$ bi--modules and admit filtrations $F_\cdot$, and
$G_\cdot$ by total order (resp. horizontal order);
in particular $G_0(\Cal A^{i,j})$ consists of the $S_j$--linear
or `relative' operators on $\Cal E_j$ and
$G_0(\Cal A^{i,j})/\goth mG_0A^{i,j} = \Cal U_{i-j}(\g)$.
Put $\bar\Cal A^{i,j}=\Cal A^{i,j}/\goth mG_0(A^{i,j})$.  As in [R1] define
$$
W_m^{i,-j} = \lambda^{j-i-1}(\g)\boxtimes
F_m(\sigma^{i+1}(\bar\Cal A^{m,j-1})),
\kern2em(i,j)\neq(m-2,m)
$$
with $\sigma^\cdot=$ exterior symmetric product (over $\Bbb C$), and
$$
W_m^{m-2,-m} = \ker(\Cal A^{m,m-1}\boxtimes F_m
(\sigma^{m-1}(\bar\Cal A^{m,m-1}))\rightarrow \sigma^m
(B^{m,m-1})\rightarrow B^m).
$$

The let $W^{\cdot,\cdot}_{m,0}$ be the analogous complex with
$\Cal A^{i,j}$ replaced by $\Cal A_0^{i,j}$ and put
$$
K^{i,-j}_m = W^{i,-j}_{m,0}\boxtimes_{S_m}\Cal E_m,\kern5em(i,j)\neq(0,0)
$$
$$
K^{0,0}_m = E.
$$

As in [R1], $W^{\cdot,\cdot}_m$ fits in an exact sequence as in (2.1),
and using the fact that $\Cal E_m$ is $S_m$--flat, $K^{\cdot,\cdot}_m$
fits in (2.2).

For instance, (2.2) has the form
$$
\CD
0 @>>> \g\boxtimes E @>>> A^{1,0}_0\boxtimes_{S_1}\Cal E @>>>
B^1(\Cal E)@>>> 0 \\
&& @VVV @VVV \\
&& E @= E \\
&& @| @| \\
&& J^\cdot_1(\g,E) && K^{\cdot,\cdot}_1 \\
\endCD
$$

\medskip

\subheading{3.  Closed 2--forms}

\medskip

Let $M$ be a manifold parameterizing a deformation of a vector bundle $E$
on $X$ with $H^0(\goth{sl}(E))=0$. For a point $0\in M$ we have the
Kodaira--Spencer map
$$
\alpha_1: T_0M\rightarrow H^1(\g), \kern2em \g=\goth{sl}(E)
$$
and its higher order analogues.  For simplicity we assume $E$ is
unobstructed, which means that the spectral sequence
$$
E_1^{p,q}=\Bbb H^q(J_m^p(\g))\Rightarrow
\Bbb H^i(J^\cdot_m(\g))
$$
degenerates at $E_1$ for $p+q=0$ for all $m$; in particular the map
$S^2H^1(\g) \rightarrow H^2(\g)$
induced by bracket vanishes.  Now we have a trace or Killing form
$$
\tr: \g\boxtimes\g\rightarrow \O_X
$$
$$
\tr(A,B) = \tr(AB).
$$
As is well known, this is a non--degenerate symmetric bilinear form and
induces an alternating form
$$
\tau:\bigwedge^2H^1(\g)\rightarrow H^2(\O_X)
$$
hence a map
$$
\tau_{M,0}:\bigwedge^2T_{M,0}\rightarrow H^2(\O_X)
$$
which we view as yielding a 2--form $\tau_M$ on $M$ with values
in the fixed vector space $H^2(\O_X)$.  When $n=\text{\rm dim}X=2$
and $K_X = O_X$, the isomorphism $\g\rightarrow \g^*$
induced by $\tr$ yields an isomorphism
$H^1(\g) \rightarrow H^1(\g^*) = H^1(\g)^*$, hence
the 2--form $\tau$ is non--degenerate on the moduli space.  In
general, any $\eta \in H^{n,n-2}(X)$ yields via cup--product with
$\eta$, a (scalar--valued) 2--form on $M$.

\proclaim{Theorem 3.1}
$\tau_M$ is a closed 2--form on $M$.
\endproclaim

{\it proof.}  It suffices to prove this for the universal second--order
deformation.  $\tr$ yields a map
$$
\tr: \sigma^2\g\rightarrow \O_X
$$
which on cohomology gives rise to
$$
\tau: H^2(\sigma^2\g) = \bigwedge^2H^1(\g)\rightarrow H^2(\O_X).
$$

Now we may identify $B^1(\bigwedge^2T_M)_0$ as
$\Bbb H^0(J_1(\g,\sigma^2\g))$, i.e. as $\Bbb H^0$
of the following complex on $\X2\times X$:
$$
\g\boxtimes\sigma^2\g\rightarrow \sigma^2\g
$$
$$
A\times(B,C) \mapsto ([A,B],C)+(C,[A,B])+(B,[A,C])+([A,C],B)
$$
Now applying $\tr$ we get
$$
2\tr(ABC-BAC+BAC-BCA)=2\tr(A(BC)-(BC)A)=0,
$$
hence $\tr$ is a map of $\g$--modules, where $\g$ acts trivially
on $\O_X$, hence $\tr$ extends to a map of complexes
$$
J^\cdot_1(\g,\sigma^2\g)\rightarrow J^\cdot_1(\g,\O).
$$
Now the differential $\g\boxtimes\sigma^2\g
\rightarrow\sigma^2\g$ clearly vanishes on $\sigma^3\g$,
so that $\sigma^3\g[2]$ forms a sub--complex of
$J_1(\g,\sigma^2\g)$ and
we have a diagram
$$
\CD
\sigma^3(\g) \\
@VVV \\
\g\boxtimes\sigma^2\g @>>> \sigma^2(\g) \\
@VVV @VVV \\
\lambda^2(\g)\boxtimes\g @>>> \g\boxtimes\g \\
\endCD
$$
where the left vertical arrows compose to zero and $\Bbb H^0$ of the
bottom complex is $T^2_{R_2}\otimes_{R_2}T^1$; in fact
the bottom complex splits as a direct sum of
$J^\cdot_3(\g)/J^\cdot_1(\g)$ and a complex
$C^\cdot:\sigma^{2,1}(\g)\rightarrow\sigma^2(\g)$,
$\sigma^{2,1}$ being the mixed tensor power
and $J^1_\cdot(\g,\sigma^2(\g))=\sigma^3(\g)[2]\oplus C^\cdot$.
As discussed in $\S$1.3,
$d\tau: \bigwedge^3H^1(\g)\rightarrow H^2(\O_X)$ coincides with the
composite
$$
\multline
\bigwedge^3H^1(\g)=H^3(\sigma^3\g)\rightarrow
\Bbb H^2(J_1(\g,\sigma^2\g)\rightarrow
\\
\Bbb H^2(\g,\O) = H^2(\O)\oplus H^1(\g)\otimes H^2(\O)\rightarrow
H^2(\O)
\endmultline
$$
But in view of the diagram
$$
\CD
\sigma^3\g[2] @>>> J_1(\g,\sigma^2\g) \\
@VVV @VVV \\
\g\boxtimes\O[2] @>>> J^\cdot_1(\g,\O) @>>> \O[1] \\
\endCD
$$
and the fact that the bottom arrows compose to zero,
$d\tau$ clearly vanishes. \qed

\medskip

\remark{Remark}
If $X$ admits a symplectic form $\omega$, and is moreover
K\"ahler, it is shown by Mukai [M] for $n=2$, and Kobayashi [K]
in general that the 2--form $\tau\wedge\omega^{n/2-1}\wedge
\bar\omega^{n/2}$ is non--degenerate (and closed), yielding
a symplectic structure on the moduli space.
\endremark

\medskip

\subheading{References}
\roster
\medskip

\item"[K]"  Kobayashi, S.: `Differential Geometry of Complex Vector
Bundles'.  I\-wa\-na\-mi and Princeton University Press, 1987.

\item"[M]"  Mukai, S.:  `Symplectic Structure on the Moduli Space of
Sheaves on an Abelian or K3 Surface'.  Invent.~Math.
{\bf 77} (1984), 101-116.

\item"[R1]"  Ran, Z.:  `Canonical Infinitesimal Deformations'.
preprint.

\item"[R2]"  $\underline{\text{\hskip4em}}$:
`Infinitesimal Deformations of Vector Bundles and Their Cohomology
Groups'.  (in preparation).

\endroster

\enddocument